\documentclass[reqno,12pt]{amsart}
\usepackage{amsmath, amssymb, amsthm, amsfonts, mathtools} 
\usepackage[english]{babel}
\usepackage{bbm}
\usepackage{graphicx}
\usepackage{url}
\usepackage{soul}
\usepackage{epstopdf}
\usepackage[ruled,vlined]{algorithm2e}
\usepackage[a4paper,bindingoffset=0.5cm,left=2cm,right=2cm,top=2.5cm,bottom=2cm,footskip=.8cm]{geometry}
\usepackage{rotating}
\usepackage{amsbsy,enumerate}
\usepackage{comment}
\usepackage{mathrsfs} 
\usepackage{enumitem}
\usepackage{cancel}
\usepackage{subcaption}
\usepackage{mathrsfs}
\usepackage{float}
\usepackage{xcolor} 
\usepackage{dashrule}

\newcommand{\bs}{\boldsymbol}

\newcommand{\vb}{\vspace{3.2mm}}
\renewcommand{\hat}{\widehat}

\newcommand{\vertiii}[1]{{\left\vert\kern-0.25ex\left\vert\kern-0.25ex\left\vert #1 \right\vert\kern-0.25ex\right\vert\kern-0.25ex\right\vert}}
\allowdisplaybreaks
\usepackage{hyperref}

\allowdisplaybreaks

\newtheorem{lemma}{Lemma}

\newtheorem{theorem}{Theorem}
\newtheorem{remark}{Remark}
\newtheorem{example}{Example}

\usepackage[utf8]{inputenc}
\usepackage{type1cm}         
\usepackage{multicol}        
\usepackage[bottom]{footmisc}

\usepackage{newtxtext}       %
\usepackage{newtxmath}       

\usepackage{pgfplots}

\usepackage{tikz}
\usetikzlibrary{automata,positioning}
\usetikzlibrary{decorations.pathreplacing,decorations.markings,shapes.geometric}
\usetikzlibrary{shapes,arrows}
\usetikzlibrary{backgrounds,calc,positioning}

\usepackage{pgfplots}
\pgfplotsset{compat=1.9}

\usetikzlibrary{chains,shapes.multipart}
\usetikzlibrary{shapes,calc}
\usetikzlibrary{automata,positioning}

\usepackage{tikz}
\usetikzlibrary{automata,positioning}
\usetikzlibrary{decorations.pathreplacing,decorations.markings,shapes.geometric}
\usetikzlibrary{shapes,arrows}
\usetikzlibrary{backgrounds,calc,positioning}

\newcommand{\twoway}[2]{%
    \draw[<->, lambda,bend left=6] (#1) to (#2);
    \draw[<->, reply,bend right=6.5] (#1) to (#2);
}

\newcommand{\oneway}[2]{%
    \draw[->, lambda, brown, bend left=7] (#1) to (#2);
    \draw[<->, reply,bend right=7] (#1) to (#2);
}

\newcommand{\fivenodes}{
\node (1) at (90:2.5) {$1$};
\node (2) at (18:2.5) {$2$};
\node (3) at (-54:2.5) {$3$};
\node (4) at (-126:2.5) {$4$};
\node (5) at (162:2.5) {$5$};
}

\tikzset{
    lambda/.style={
        thick,
        shorten >=2pt,
        shorten <=2pt
    },
    reply/.style={
        thick,
        dashed,
        shorten >=2pt,
        shorten <=2pt
    }
}

\tikzset{
  every node/.style={
    circle,
    draw,
    fill=blue!30,
    minimum size=6mm
  }
}
\usepackage{xcolor}
\usepackage{tikz}

\definecolor{plotred}{RGB}{255,0,0}
\definecolor{plotblue}{RGB}{0,0,255}
\definecolor{plotgreen}{RGB}{0,255,0}
\definecolor{plotorange}{RGB}{255,165,0}
\definecolor{plotpurple}{RGB}{160,32,240}

\newcommand{\Rred}{\tikz[baseline=-0.6ex]{
\draw[plotred,line width=0.8pt] (0,0)--(0.5,0);
\fill[plotred] (0.25,0) circle (2.5pt);
}}

\newcommand{\Rblue}{\tikz[baseline=-0.6ex]{
\draw[plotblue,line width=0.8pt] (0,0)--(0.5,0);
\fill[plotblue] (0.25,0.06)--(0.19,-0.05)--(0.31,-0.05)--cycle;
}}

\newcommand{\Rgreen}{\tikz[baseline=-0.6ex]{
\draw[plotgreen,line width=0.8pt] (0,0)--(0.5,0);
\fill[plotgreen] (0.20,-0.05) rectangle (0.30,0.05);
}}

\newcommand{\Rorange}{\tikz[baseline=-0.6ex]{
\draw[plotorange,line width=0.8pt] (0,0)--(0.5,0);
\fill[plotorange] (0.25,0.06)--(0.31,0)--(0.25,-0.06)--(0.19,0)--cycle;
}}

\newcommand{\Rpurple}{\tikz[baseline=-0.6ex]{
\draw[plotpurple,line width=0.8pt] (0,0)--(0.5,0);
\draw[plotpurple,line width=0.8pt]
(0.18,0)--(0.32,0)
(0.25,-0.07)--(0.25,0.07)
(0.20,-0.05)--(0.30,0.05)
(0.20,0.05)--(0.30,-0.05);
}}

\begin{document}

	\title[Inference in an interactive multi-type branching network]{Inference in an interactive multi-type\\ branching network with partial information}
\author[G. Hasanli, M. Mandjes, {\tiny and} J. Wang]{G. Hasanli, M. Mandjes, {\tiny and} J. Wang}


\begin{abstract}
	We study statistical inference for a high-dimensional stochastic interaction model observed only through low-dimensional aggregate data. The model describes pairwise communication between $d$ individuals as a collection of mutually interacting branching processes with immigration, where activity on one directed pair may trigger subsequent activity in the reverse direction. The latent interaction counts are not observed. Instead, for each individual, only the total number of messages sent and received in each time slot is available. This aggregation gives rise to an inverse problem in which the parameters governing the latent pairwise dynamics must be recovered from substantially compressed observations.

\noindent We establish identifiability of the model parameters under a tractable parametric specification and develop a method-of-moments estimator based on stationary means and lagged covariances of the observed aggregate processes. In the leading Poisson--Bernoulli specification, the relevant stationary moments admit explicit expressions, which allow the latent interaction parameters to be recovered from a finite collection of observable moments. We prove consistency and asymptotic normality of the resulting estimator. The estimation procedure avoids likelihood maximization over the latent high-dimensional process and remains computationally simple despite the loss of information caused by aggregation.

\noindent Although motivated by communication networks, the framework applies more broadly to partially observed interacting count processes in which fine-scale dynamics are latent and only aggregate measurements are available.

\vb


\vb

\noindent
{\sc Affiliations.} 
GH is  with Mathematical Institute, Leiden University, Leiden, 
The Netherlands.

\noindent 
MM is  with Mathematical Institute, Leiden University, Leiden, 
The Netherlands; he is also affiliated with (a)~Korteweg-de Vries Institute for Mathematics, University of Amsterdam, Amsterdam, The Netherlands, (b)~E{\sc urandom}, Eindhoven University of Technology, Eindhoven, The Netherlands, (c)~Amsterdam Business School, Faculty of Economics and Business, University of Amsterdam, Amsterdam, The Netherlands.  

\noindent JW is with School of Engineering, Macquarie University, Sydney, Australia.

\noindent Date: {\it \today}.

\vb

\noindent
{\sc Acknowledgments.} 
MM's and JW's research has been funded by the NWO Gravitation project {\tiny NETWORKS} under grant agreement no.\ 024.002.003.

\vb

\noindent
{\sc Email.} \url{g.s.hasanli@umail.leidenuniv.nl}, \:\url{m.r.h.mandjes@math.leidenuniv.nl}, 

\noindent \url{jiesenwang@gmail.com}

\end{abstract}

 \maketitle

\newpage
\section{Introduction}\label{sec:intro}

Interactions between individuals in digital communication systems --- such
as email exchanges, messaging platforms, and online discussion
threads --- generate complex dependence structures that evolve over time.
Activity between two individuals may trigger subsequent activity in the
opposite direction, while interaction intensities may differ across
pairs. Statistical inference becomes particularly challenging when the
individual pairwise interactions are not directly observed.

A natural question is whether the parameters governing such latent
interactions can be recovered from aggregate activity. This problem
arises in communication networks, but also in network traffic monitoring,
partially logged communication systems, and other interacting count
processes in which observations are available only after aggregation.
The central statistical difficulty is that aggregation substantially
reduces the dimension of the available data while retaining dependence
generated by the unobserved interaction mechanisms.

\smallskip

{\it Model \& objectives ---}
In this paper, we study a stochastic model for pairwise interactions
among $d\geqslant3$ individuals. For each ordered pair $(i,j)$ with
$i\neq j$, the latent process records the number of interactions from
$i$ to $j$ at each time point. Thus, the latent system contains
$d^2-d$ pairwise processes. We assume, however, that these processes are
not observed individually. Instead, for each individual we observe only
the total outgoing and incoming activity. The observation dimension is therefore $2d$, compared with the latent
dimension $d^2-d$. The two dimensions coincide when $d=3$, while for
$d>3$ the observation scheme strictly reduces the dimension of the
available data.

The main question of the paper is whether the parameters governing the
latent pairwise dynamics can nevertheless be identified and estimated
from these aggregate observations. Our approach is based on stationary
moment relations. In particular, while marginal moments generally
contain only aggregate information, lagged covariances retain information
about the direction and strength of the interaction mechanism. We exploit
these temporal relations to construct an explicit identification map
and a corresponding method-of-moments estimator.

We focus on a parametric specification in which the immigration and
offspring distributions are each indexed by a single parameter. The
main example is Poisson immigration combined with Bernoulli offspring.
Under the corresponding stability conditions, we show that a finite
collection of observable moments identifies all latent parameters. The
estimator is obtained directly from empirical moments of the aggregate
observations and therefore does not require reconstruction of the latent
pairwise processes or maximization of a high-dimensional latent-data
likelihood.

The resulting procedure has two main advantages. First, it uses only the
aggregate observations that are available in the proposed observation
scheme. Second, the parameter recovery is explicit, so that the
estimator can be implemented by evaluating empirical moments and applying
the inverse moment map. We establish consistency and asymptotic
normality under the Poisson--Bernoulli specification.

\smallskip

{\it Literature \& contribution ---}
The model considered in this paper belongs to the broad class of
branching processes with immigration. Classical foundations for
branching processes with immigration can be found in
\cite{AthreyaNey1972,Mode}, while multitype branching processes and their
long-term behavior are treated in detail in
\cite{Bass,Harris1963}. Related interaction structures have been studied
in coupled branching systems and multitype population models
\cite{Biggins,Jagers,Nerman}.

A second relevant strand of literature concerns statistical inference
under partial or indirect observation. In many stochastic systems,
individual-level processes are latent while only aggregated measurements
are available. Such problems arise, for example, in ecological models,
epidemiology, queueing networks, and other partially observed
stochastic systems \cite{Diggle,GuttorpThompson,RavnerBoxmaMandjes,vanKampen}.
Related aggregation problems also arise in dynamic network models; see,
for example, \cite{MandjesWangERGM}. The present setting differs in that
the aggregation is applied to a high-dimensional interacting branching
system, while the objective is to recover parameters governing the
underlying pairwise dynamics.

From a statistical perspective, inference for branching systems has been
studied using likelihood-based methods, estimating equations, and
moment-based procedures. Moment methods are particularly natural when
the latent process is high-dimensional or only indirectly observed. An
example of inference for partially observed multitype branching systems
is given in \cite{Laredo}.

Our contribution is to study parameter identification and estimation for
a fully coupled $(d^2-d)$-dimensional interaction system when only $2d$
aggregate processes are observed. We derive observable moment relations
that retain information about the latent pairwise interactions,
establish identification under the Poisson--Bernoulli specification, and
construct an explicit method-of-moments estimator. We further establish
consistency and asymptotic normality of the resulting estimator.

The identification argument also illustrates a broader point:
aggregation does not necessarily eliminate information about latent
interaction mechanisms. In the present model, temporal dependence
generated by the alternating interaction between the two directions of
each pair remains visible in the aggregate observations and can be used
for parameter recovery.

\smallskip

{\it Organization ---}
Section \ref{sec: model} introduces the latent interaction model, the
observation scheme, and the standing assumptions. Section
\ref{sec: mom} derives the stationary moment relations used for
identification, including the lag-$0$, lag-$1$, and lag-$2$ moments of
the observed processes. Section \ref{sec: as} develops the
method-of-moments estimator and establishes identification, consistency,
and asymptotic normality. Section \ref{sec:num} presents numerical
experiments illustrating the finite-sample behavior of the estimator.

\section{Model and notations}\label{sec: model}

We consider a system of $d\in\mathbb N$ individuals who interact by
exchanging messages. For each ordered pair $(i,j)$ with $i\neq j$, let
$X_{ij,n}$ denote the number of distinct message threads on which a
message is sent from individual $i$ to individual $j$ during time slot
$n$. The latent system is therefore described by the
$(d^2-d)$-dimensional process
\[
(X_{ij,n})_{i\neq j,\;n\in\mathbb Z}.
\]

In every time slot, new threads may start spontaneously, while existing
threads may generate subsequent messages in the opposite direction. We
model these two mechanisms through the branching process
\begin{equation}\label{eq:model}
X_{ij,n+1}
=
I_{ij,n+1}
+
\sum_{k=1}^{X_{ji,n}}\xi_{ji,n+1,k},
\qquad i\neq j.
\end{equation}
Here, $I_{ij,n}$ represents exogenous immigration into the directed pair
$(i,j)$, while $\xi_{ij,n,k}$ represents the number of offspring
generated by the $k$th interaction from $i$ to $j$. Thus, activity from
$j$ to $i$ may generate subsequent activity from $i$ to $j$, and
conversely.

The immigration variables $(I_{ij,n})_{n\in\mathbb Z}$ are assumed to be
independent across $(i,j,n)$ and identically distributed over time, with
non-negative integer-valued generic variable $I_{ij}$ and finite mean
\[
\lambda_{ij}:=\mathbb E[I_{ij}].
\]
The offspring variables $(\xi_{ij,n,k})$ are likewise assumed
independent across $(i,j,n,k)$ and independent of the immigration
variables. They are identically distributed over time and across
individuals of a given ordered pair, with non-negative integer-valued
generic variable $\xi_{ij}$ and finite mean
\[
m_{ij}:=\mathbb E[\xi_{ij}].
\]

The latent system decomposes naturally over unordered pairs. For each
unordered pair $\{i,j\}$, the bivariate process
$(X_{ij,n},X_{ji,n})_{n\in\mathbb Z}$ evolves independently of the
processes associated with all other unordered pairs. Over two time
steps, the corresponding mean reproduction factor is
$m_{ij}m_{ji}$. We therefore impose the stability condition
\[
m_{ij}m_{ji}<1,
\qquad i\neq j,
\]
and work throughout with the stationary version of the process.

The immigration and offspring distributions are assumed to belong to
parametric families. We focus on the case in which each family is
indexed by a single parameter, so that the higher moments of $I_{ij}$
and $\xi_{ij}$ are determined by $\lambda_{ij}$ and $m_{ij}$,
respectively. This restriction keeps the moment equations and
identification argument explicit. The resulting quantities
\[
\{\lambda_{ij},m_{ij}:i\neq j\}
\]
thus constitute the $2d^2-2d$ unknown parameters of interest.

\begin{example}\label{ex1}
{\em
The main specification considered in this paper combines Poisson
immigration with Bernoulli offspring:
\[
I_{ij}\sim\operatorname{Poisson}(\lambda_{ij}),
\qquad
\xi_{ij}\sim\operatorname{Bernoulli}(p_{ij}),
\]
so that $m_{ij}=p_{ij}$. As we will see, under this specification, the expressions for the stationary
moments simplify substantially, allowing the latent parameters to be
recovered explicitly from a finite collection of observable moments.
}
\hfill$\Diamond$
\end{example}

Importantly, the latent processes $X_{ij,n}$ are {\it not} observed.
Instead, the data consist only of the total outgoing and incoming
activity of each individual. In other words, for each $i=1,\ldots,d$, we
observe the two processes
\begin{equation}\label{obs}
Y_{i,n}:=\sum_{j\neq i}X_{ij,n},
\qquad
Z_{i,n}:=\sum_{j\neq i}X_{ji,n}.
\end{equation}
Thus, $Y_{i,n}$ records the total number of threads sent by individual
$i$, while $Z_{i,n}$ records the total number of threads received by
individual $i$ during time slot $n$. These $2d$ aggregate processes are
the only data available for inference; the individual pairwise counts
$X_{ij,n}$ are latent.

The observation scheme thus replaces the $d^2-d$ latent pairwise processes by
only $2d$ aggregate processes. The statistical objective is therefore
to recover the parameters $\lambda_{ij}$ and $m_{ij}$ from the observed
processes in \eqref{obs}.

\begin{remark}\label{R1}
{\em
Before aggregation, each unordered pair $\{i,j\}$ forms a bivariate
Markov chain, and different unordered pairs are mutually independent.
After aggregation, these properties are no longer directly available:
each observed component combines several latent pairwise processes, and
the resulting aggregate process is generally not Markovian. In
particular, likelihood-based procedures based directly on the latent
transition structure cannot be applied to the aggregate observations
without accounting for the unobserved pairwise processes.

This motivates the moment-based approach developed below. The stationary
temporal moments of the aggregate observations retain information about
the underlying pairwise interactions and will be used to identify and
estimate the latent parameters.
}
\hfill$\Diamond$
\end{remark}

The next section derives the stationary moments of the latent pairwise
processes and expresses the corresponding moments of the aggregate
observations in terms of the unknown parameters. These relations form
the basis of the identification and estimation procedures developed
subsequently.
\section{Moments in stationarity}\label{sec: mom}

Our estimation procedure is based on stationary moments of the aggregate observations. We therefore first derive the relevant moment relations for the latent pairwise processes and then determine which of these relations are available from the aggregated observations. Throughout this section, the process is assumed to be in its stationary regime. Throughout, we let $X_{ij}$ denote a generic random variable with the stationary distribution of $X_{ij,n}$.

\subsection{Lag-0 moments}\label{ssec: lag0}

We begin by computing the stationary means. These follow from the standard
mean-balance equations for a branching process with immigration. Write
$\mu_{ij}:=\mathbb E[X_{ij}]$. 
Applying Wald's equation to \eqref{eq:model} gives
\begin{equation}\label{means}
\mu_{ij}=\lambda_{ij}+\mu_{ji}m_{ji},
\qquad
\mu_{ji}=\lambda_{ji}+\mu_{ij}m_{ij}.
\end{equation}
Solving this system yields
\begin{equation}\label{means-sol}
\mu_{ij}
=
\frac{\lambda_{ij}+\lambda_{ji}m_{ji}}
{1-m_{ij}m_{ji}},
\qquad
\mu_{ji}
=
\frac{\lambda_{ji}+\lambda_{ij}m_{ij}}
{1-m_{ij}m_{ji}}.
\end{equation}
Recall that $m_{ij}m_{ji}<1$ is the stability condition for the unordered pair $\{i,j\}$.

The stationary means of the observed processes are obtained by summing the latent means. With $\mu_i^Y:=\mathbb E[Y_i]$ and $\mu_i^Z:=\mathbb E[Z_i]$, we have \[\mu_i^Y=\sum_{j\neq i}\mu_{ij},\qquad \mu_i^Z=\sum_{j\neq i}\mu_{ji}.\] Thus, the observed first moments provide $2d$ equations for the latent means.

\begin{remark}\label{R2}
{\em
The $2d$ first-moment equations are not all independent. Indeed, at every time point the total outgoing activity equals the total incoming activity,
\[
\sum_{i=1}^dY_i=\sum_{i=1}^dZ_i.
\]
Consequently, $\sum_{i=1}^d\mu_i^Y=\sum_{i=1}^d\mu_i^Z$, and the first moments provide at most $\ell_0:=2d-1$ independent equations.
}
\hfill$\Diamond$
\end{remark}

\begin{remark}\label{R3}
{\em
Under the Poisson--Bernoulli specification of Example~\ref{ex1}, the stationary variables $X_{ij}$ and $X_{ji}$ are Poisson distributed. Indeed, the stationary probability generating function equations imply
\[
X_{ij}\sim\operatorname{Poisson}(\mu_{ij}),
\qquad
X_{ji}\sim\operatorname{Poisson}(\mu_{ji}),
\]
with means given by \eqref{means-sol}. In particular, all moments of the stationary latent variables are finite.
}
\hfill$\Diamond$
\end{remark}

Moments pertaining only to the marginal stationary distributions cannot
recover the individual pairwise parameters. This is particularly
transparent under the Poisson--Bernoulli specification. By
Remark~\ref{R3}, each latent variable is Poisson, and independence across
distinct unordered pairs implies that the aggregate variables $Y_i$ and
$Z_i$ are also Poisson, with means $\mu_i^Y$ and $\mu_i^Z$. Consequently,
all stationary marginal moments of $Y_i$ and $Z_i$ are determined by their
means. Thus, estimating (say) stationary second moments or higher marginal
moments does not provide information beyond the first-moment equations
about the individual $\lambda_{ij}$ and $m_{ij}$. We therefore turn to
temporal dependence. The mutual excitation mechanism leaves information
in lagged covariances even though the latent pairwise counts are
unobserved.

\subsection{Lag-1 moments}\label{ssec:lag1}

We first consider lag-1 covariances of the aggregate processes. By stationarity, these quantities do not depend on the time index. We write 
\[\rho_{ij}^Y:=\operatorname{Cov}(Y_{i,n},Y_{j,n+1}),\qquad \rho_{ij}^Z:=\operatorname{Cov}(Z_{i,n},Z_{j,n+1}).\]
The covariance relations are most conveniently derived at the latent
pair level. For a fixed unordered pair $\{i,j\}$, let
$V_n:=(X_{ij,n},X_{ji,n})^\top$ and define
$\Sigma_k:=\operatorname{Cov}(V_n,V_{n+k})$. By \eqref{eq:model},
the conditional mean of $V_{n+1}$ given $V_n$ is
\[
\mathbb E[V_{n+1}\mid V_n]
=
\begin{pmatrix}
\lambda_{ij}\\
\lambda_{ji}
\end{pmatrix}
+
A V_n,
\qquad
A:=
\begin{pmatrix}
0&m_{ji}\\
m_{ij}&0
\end{pmatrix}.
\]
The alternating interaction structure is therefore encoded in the
matrix $A$. The following covariance recursion is then an immediate
consequence of the conditional mean relation above.

\begin{lemma}\label{L2}
Let $V_n=(X_{ij,n},X_{ji,n})^\top$ and define $\Sigma_k:=\operatorname{Cov}(V_n,V_{n+k})$. Then, for every $k\geqslant0$,
\[
\Sigma_k=\Sigma_0A^k,
\qquad
A=
\begin{pmatrix}
0&m_{ji}\\
m_{ij}&0
\end{pmatrix}.
\]
\end{lemma}

\begin{proof}
From the conditional mean relation above,
\[
\Sigma_1
=
\operatorname{Cov}(V_n,V_{n+1})
=
\Sigma_0A.
\]
Iterating the same argument gives $\Sigma_k=\Sigma_0A^k$ for every $k\geqslant1$.
\end{proof}

At a common time point, the two components of a latent pair are independent; see Lemma \ref{L3} in the appendix. As a consequence,  $\Sigma_0$ is diagonal. Write $\sigma_{ij}^2:=\operatorname{Var}(X_{ij})$ and $\sigma_{ji}^2:=\operatorname{Var}(X_{ji})$. To obtain these variances, define
\[
\phi_{ij}
:=
\mathbb E[I_{ij}^2]
+2\lambda_{ij}\mu_{ji}m_{ji}
+\mu_{ji}\operatorname{Var}(\xi_{ji}).
\]
Taking squares in \eqref{eq:model}, using independence of the innovation terms from the current state, and imposing stationarity gives $\mathbb E[X_{ij}^2]=\phi_{ij}+m_{ij}^2\mathbb E[X_{ji}^2]$ and the analogous equation with $i$ and $j$ interchanged. Solving the resulting system yields
\begin{equation}\label{exij2}
\mathbb E[X_{ij}^2]
=
\frac{\phi_{ij}+\phi_{ji}m_{ji}^2}
{1-m_{ij}^2m_{ji}^2},
\qquad
\mathbb E[X_{ji}^2]
=
\frac{\phi_{ji}+\phi_{ij}m_{ij}^2}
{1-m_{ij}^2m_{ji}^2}.
\end{equation}
We have thus identified $\sigma_{ij}^2=\mathbb E[X_{ij}^2]-\mu_{ij}^2$, and similarly $\sigma_{ji}^2=\mathbb E[X_{ji}^2]-\mu_{ji}^2$.

\begin{remark}\label{R5}
{\em
Since $A$ has zero diagonal entries, the covariance matrices alternate between diagonal and zero-diagonal structures: $\Sigma_k$ is diagonal for even $k$ and has zero diagonal for odd $k$. Moreover, the spectral radius of $A$ is $\sqrt{m_{ij}m_{ji}}<1$ under the stability condition. Hence the latent covariance sequence decays geometrically.
}
\hfill$\Diamond$
\end{remark}

We can now express the lag-1 covariances of the aggregate observations in terms of the latent pairwise covariances. For $i\neq j$, independence between distinct unordered pairs implies that only the pair $\{i,j\}$ contributes to $\rho_{ij}^Y$. Specifically, $\rho_{ij}^Y=\operatorname{Cov}(X_{ij,n},X_{ji,n+1})$, which is the $(1,2)$ entry of $\Sigma_1$. Since $\Sigma_0$ is diagonal, this gives $\rho_{ij}^Y=m_{ji}\sigma_{ij}^2$. For $i=j$, the corresponding entry of $\Sigma_1$ is zero, so $\rho_{ii}^Y=0$. Substituting the stationary mean and second-moment expressions gives the following explicit relation.

\begin{lemma}\label{L5}
Fix $i,j\in\{1,\ldots,d\}$. If $i\neq j$ and
$m_{ij}m_{ji}<1$, then
\[
\rho_{ij}^Y
=
\mu_{ij}\lambda_{ji}
+
m_{ji}
\frac{\phi_{ij}+\phi_{ji}m_{ji}^2}
{1-m_{ij}^2m_{ji}^2}
-
\mu_{ij}\mu_{ji}.
\]
Moreover, $\rho_{ii}^Y=0$.
\end{lemma}

The corresponding expression for $\rho_{ij}^Z$ is obtained by interchanging $i$ and $j$. Thus $\rho_{ij}^Z=\rho_{ji}^Y$, so the lag-1 covariances of $Z$ contain no additional information beyond those of $Y$.

\begin{remark}\label{R7}
{\em
The relation $\rho_{ij}^Z=\rho_{ji}^Y$ implies that the lag-1 covariances of $Y$ and $Z$ provide at most $\ell_1^\leftrightarrow:=d^2-d$ independent equations.
}
\hfill$\Diamond$
\end{remark}

For completeness, consider also the mixed lag-1 covariances \[\rho_{ij}^{YZ}:=\operatorname{Cov}(Y_{i,n},Z_{j,n+1}), \qquad \rho_{ij}^{ZY}:=\operatorname{Cov}(Z_{i,n},Y_{j,n+1}).\] For $i\neq j$, the only potentially contributing latent covariance corresponds to a diagonal entry of $\Sigma_1$ and therefore vanishes. Hence $\smash{\rho_{ij}^{YZ}=\rho_{ij}^{ZY}=0}$ for $i\neq j$. For $i=j$, the mixed covariances reduce to sums of the lag-1 covariances already considered: $\smash{\rho_{ii}^{YZ}=\sum_{k\neq i}\rho_{ik}^Y}$ and $\smash{\rho_{ii}^{ZY}=\sum_{k\neq i}\rho_{ik}^Z}$.

\begin{remark}\label{R8}
{\em
Thus, the mixed lag-1 covariances provide no additional independent moment equations.
}
\hfill$\Diamond$
\end{remark}

\subsection{Lag-2 moments}

Upon combining Remarks \ref{R2}, \ref{R7} and \ref{R8}, we conclude that the lag-0 and lag-1 moments  provide at most $\ell_0+\ell_1^\leftrightarrow=(2d-1)+(d^2-d)=d^2+d-1$ independent equations. Since there are $2d^2-2d$ unknown parameters, these moments alone cannot identify the parameters. We therefore turn to lag-2 covariances.

Define \[\rho_{ij,2}^Y:=\operatorname{Cov}(Y_{i,n},Y_{j,n+2}),\qquad\rho_{ij,2}^Z:=\operatorname{Cov}(Z_{i,n},Z_{j,n+2}).\] Since $\Sigma_2=\Sigma_0A^2$ is diagonal, these covariances vanish for $i\neq j$. For $i=j$, the lag-2 covariance of $Y_i$ is obtained by summing the corresponding latent covariances, giving $\rho_{ii,2}^Y=\sum_{k\neq i}m_{ik}m_{ki}\sigma_{ik}^2$. The corresponding expression for $Z_i$ is obtained by relabelling the indices.

The useful additional information comes from the mixed lag-2 covariances. Define \[\rho_{ij,2}^{YZ}:=\operatorname{Cov}(Y_{i,n},Z_{j,n+2}),\qquad \rho_{ij,2}^{ZY}:=\operatorname{Cov}(Z_{i,n},Y_{j,n+2}).\] For $i=j$, the relevant latent covariance is an off-diagonal entry of $\Sigma_2$ and therefore vanishes. For $i\neq j$, independence between distinct unordered pairs leaves only the component $X_{ij}$, so
\[
\rho_{ij,2}^{YZ}
=
\operatorname{Cov}(X_{ij,n},X_{ij,n+2})
=
m_{ij}m_{ji}\sigma_{ij}^2.
\]
Thus the lag-2 covariance returns to the same directed component and is multiplied by the product $m_{ij}m_{ji}$. Substituting \eqref{means} and \eqref{exij2} gives the explicit expression below.

\begin{lemma}
Fix $i,j\in\{1,\ldots,d\}$. If $i\neq j$ and
$m_{ij}m_{ji}<1$, then
\[\rho_{ij,2}^{YZ}=\left(\lambda_{ij} \mu_{ij} + m_{ji} \lambda_{ji} \mu_{ij} 
    + m_{ji} m_{ij} \frac{\phi_{ij}+\phi_{ji}m_{ji}^2}{1-m_{ij}^2m^2_{ji}}-\mu_{ij}^2\right){\bs 1}_{\{i\not=j\}}.\]
Moreover, {$\rho_{ii}^{YZ}=0$}.
\end{lemma}

The corresponding expression for $\smash{\rho_{ij,2}^{ZY}}$ is obtained
by interchanging $i$ and $j$, so
$\smash{\rho_{ij,2}^{ZY}=\rho_{ji,2}^{YZ}}$. Moreover, the diagonal
lag-2 covariances of $Y$ are already determined by the mixed
covariances, since
$\smash{\rho_{ii,2}^Y=\sum_{k\neq i}\rho_{ik,2}^{YZ}}$. Hence, the only
additional lag-2 information is provided by the $d^2-d$ mixed
covariances $\smash{\rho_{ij,2}^{YZ}}$ with $i\neq j$.

\begin{remark}\label{R9}
{\em
The lag-2 moments therefore provide at most $\ell_2^\leftrightarrow:=d^2-d$ additional independent equations. Combining the lag-0, lag-1, and lag-2 information gives, recalling Remarks \ref{R2}, \ref{R7} and \ref{R8},
\[
\ell_0+\ell_1^\leftrightarrow+\ell_2^\leftrightarrow
=
(2d-1)+2(d^2-d)
=
2d^2-1.
\]
Thus, the resulting moment system is overidentified relative to the $2d^2-2d$ unknown parameters. This count alone does not establish identification, however, since the moment equations may be functionally dependent. In the general parametric setting, identification is determined by the rank and local invertibility properties of the corresponding population moment map.
}
\hfill$\Diamond$
\end{remark}

Under the Poisson--Bernoulli specification, the moment relations simplify substantially. This allows us to select a convenient set of observable moments and derive an explicit inverse map from those moments to the four parameters associated with each unordered pair. The next section develops this inversion, establishing identification and then constructing the corresponding estimator and studying its asymptotic properties.

\section{Estimation and asymptotic normality}\label{sec: as}

We now turn to the estimation problem. Throughout this section, we consider the leading specification introduced in Example~\ref{ex1}, in which the immigration variables are Poisson and the offspring variables are Bernoulli. Thus, $I_{ij}\sim\operatorname{Poisson}(\lambda_{ij})$ and $\xi_{ij}\sim\operatorname{Bernoulli}(m_{ij})$ for $i\neq j$, and the parameter vector consists of the $2d^2-2d$ quantities $\lambda_{ij}$ and $m_{ij}$.

The discussion at the end of  Subsection~\ref{ssec: lag0} shows that marginal stationary
moments do not identify the latent parameters, so identification must
come from temporal dependence. Under the Poisson--Bernoulli
specification, the lagged covariance relations derived in
Section~\ref{sec: mom} simplify substantially because
$\operatorname{Var}(X_{ij})=\mu_{ij}$. This simplification makes it
possible to construct an explicit inverse map from a finite collection
of observable lagged moments to the latent parameters. We first
describe this inverse map and the resulting estimator, and then
establish consistency and asymptotic normality. The idea is to use the $d^2-d$ lag-1 covariances and the $d^2-d$ mixed lag-2 covariances. 
\subsection{Estimator and identifiability}
For $i\neq j$, Remark~\ref{R3} gives
$\operatorname{Var}(X_{ij})=\mu_{ij}$ under the
Poisson--Bernoulli specification. Substituting this identity into the
relations of
Section~\ref{sec: mom}, we obtain for the lag-1 covariances and the mixed lag-2 covariances the relations that
\begin{equation}\label{eq:BinPo}
\rho_{ij}^Y=\mu_{ij}m_{ij},
\qquad
\rho_{ij,2}^{YZ}=\mu_{ij}m_{ij}m_{ji},
\qquad i\neq j,
\end{equation}
respectively. 
These are the observable moment relations that we use to recover the
latent parameters, as follows.

For each unordered pair $\{i,j\}$, applying \eqref{eq:BinPo} to both
orientations gives
$\smash{\rho_{ij}^Y=\mu_{ij}m_{ij}}$,
$\smash{\rho_{ji}^Y=\mu_{ji}m_{ji}}$,
$\smash{\rho_{ij,2}^{YZ}=\mu_{ij}m_{ij}m_{ji}}$, and
$\smash{\rho_{ji,2}^{YZ}=\mu_{ji}m_{ji}m_{ij}}$.
Taking ratios therefore yields
\begin{equation}\label{eq:m_inverse}
m_{ji}
=
\frac{\rho_{ij,2}^{YZ}}{\rho_{ij}^Y},
\qquad
m_{ij}
=
\frac{\rho_{ji,2}^{YZ}}{\rho_{ji}^Y}.
\end{equation}
The stationary means are then recovered from \eqref{eq:BinPo}:
$\smash{\mu_{ij}=\rho_{ij}^Y/m_{ij}}$ and
$\smash{\mu_{ji}=\rho_{ji}^Y/m_{ji}}$. Substituting \eqref{eq:m_inverse} thus gives the relations
\[
\mu_{ij}
=
\frac{\rho_{ij}^Y\rho_{ji}^Y}{\rho_{ji,2}^{YZ}},
\qquad
\mu_{ji}
=
\frac{\rho_{ij}^Y\rho_{ji}^Y}{\rho_{ij,2}^{YZ}}.
\]
Finally, the stationary mean equations in \eqref{means} yield
$\smash{\lambda_{ij}=\mu_{ij}-\mu_{ji}m_{ji}}$ and
$\smash{\lambda_{ji}=\mu_{ji}-\mu_{ij}m_{ij}}$. As a consequence,
\begin{equation}\label{eq:lambda_inverse}
\lambda_{ij}
=
\frac{\rho_{ij}^Y\rho_{ji}^Y}
{\rho_{ji,2}^{YZ}}
-
\rho_{ji}^Y,
\qquad
\lambda_{ji}
=
\frac{\rho_{ij}^Y\rho_{ji}^Y}
{\rho_{ij,2}^{YZ}}
-
\rho_{ij}^Y.
\end{equation}
Thus, for parameters in the interior of the parameter space, the four
parameters associated with an unordered pair $\{i,j\}$ are identified
by four population moments of the observed processes, namely
$\smash{\rho_{ij}^Y}$, $\smash{\rho_{ji}^Y}$,
$\smash{\rho_{ij,2}^{YZ}}$, and
$\smash{\rho_{ji,2}^{YZ}}$. These moments are functions of the observed
processes $Y$ and $Z$ and can therefore be estimated directly from the
data.

For each unordered pair $\{i,j\}$, define
\[
{\bs\rho}_{ij}
:=
\left(
\rho_{ij}^Y,\,
\rho_{ji}^Y,\,
\rho_{ij,2}^{YZ},\,
\rho_{ji,2}^{YZ}
\right),
\qquad
{\bs\theta}_{ij}
:=
\left(
\lambda_{ij},\,
\lambda_{ji},\,
m_{ij},\,
m_{ji}
\right).
\]
Expressions \eqref{eq:m_inverse} and \eqref{eq:lambda_inverse} define an explicit inverse map
$\smash{F:\mathbb{R}^4\to\mathbb{R}^4}$, with
$\smash{{\bs\theta}_{ij}=F({\bs\rho}_{ij})}$ on the interior of the
parameter space. Equivalently, the population moment map from the four
pairwise parameters $\smash{{\bs\theta}_{ij}}$ to the four observable
moments $\smash{{\bs\rho}_{ij}}$ is one-to-one on the interior of the
parameter space, with $F$ providing its explicit inverse. Thus,
identification is pairwise: the four parameters associated with a given
unordered pair $\{i,j\}$ can be recovered from moments involving only
that pair of individuals.

We next estimate these population moments from the observed processes.
Define the empirical means
\[
\hat\mu_{i,N}^Y:=\frac{1}{N}\sum_{n=1}^N Y_{i,n},
\qquad
\hat\mu_{i,N}^Z:=\frac{1}{N}\sum_{n=1}^N Z_{i,n}.
\]
For each $i\neq j$, define the empirical cross-product averages
\[
\hat q_{ij,N}^Y
:=
\frac{1}{N-1}\sum_{n=1}^{N-1}Y_{i,n}Y_{j,n+1},
\qquad
\hat q_{ij,N}^{YZ}
:=
\frac{1}{N-2}\sum_{n=1}^{N-2}Y_{i,n}Z_{j,n+2}.
\]
The corresponding quantities with $i$ and $j$ interchanged are
$\smash{\hat q_{ji,N}^Y}$ and $\smash{\hat q_{ji,N}^{YZ}}$. The empirical
covariance estimators are therefore
\[
\hat\rho_{ij,N}^Y
=
\hat q_{ij,N}^Y-\hat\mu_{i,N}^Y\hat\mu_{j,N}^Y,
\qquad
\hat\rho_{ij,2,N}^{YZ}
=
\hat q_{ij,N}^{YZ}-\hat\mu_{i,N}^Y\hat\mu_{j,N}^Z.
\]
For each unordered pair $\{i,j\}$, define
\[
\hat{\bs\rho}_{ij,N}
:=
\left(
\hat\rho_{ij,N}^Y,\,
\hat\rho_{ji,N}^Y,\,
\hat\rho_{ij,2,N}^{YZ},\,
\hat\rho_{ji,2,N}^{YZ}
\right).
\]

For later use, we collect the empirical quantities entering these four
covariance estimators into
\[
\hat{\bs\eta}_{ij,N}
:=
\left(
\hat q_{ij,N}^Y,\,
\hat q_{ji,N}^Y,\,
\hat q_{ij,N}^{YZ},\,
\hat q_{ji,N}^{YZ},\,
\hat\mu_{i,N}^Y,\,
\hat\mu_{j,N}^Y,\,
\hat\mu_{i,N}^Z,\,
\hat\mu_{j,N}^Z
\right).
\]
For fixed $i\neq j$, let the corresponding population quantities be
\[
q_{ij}^Y:=\mathbb E[Y_{i,n}Y_{j,n+1}],
\qquad
q_{ij}^{YZ}:=\mathbb E[Y_{i,n}Z_{j,n+2}],
\qquad
\mu_i^Y:=\mathbb E[Y_{i,n}],
\qquad
\mu_i^Z:=\mathbb E[Z_{i,n}],
\]
and define
\[
\bs\eta_{ij}
:=
\left(
q_{ij}^Y,\,
q_{ji}^Y,\,
q_{ij}^{YZ},\,
q_{ji}^{YZ},\,
\mu_i^Y,\,
\mu_j^Y,\,
\mu_i^Z,\,
\mu_j^Z
\right);
\]
thus, $\bs\eta_{ij}$ is the population counterpart of
$\hat{\bs\eta}_{ij,N}$. In particular, the centering operation defining
the observable covariance vector is represented by the map
$G:\mathbb R^8\to\mathbb R^4$ given by
\[
G\left(q_{ij}^Y,q_{ji}^Y,q_{ij}^{YZ},q_{ji}^{YZ},
\mu_i^Y,\mu_j^Y,\mu_i^Z,\mu_j^Z\right)
:=
\bigl(
q_{ij}^Y-\mu_i^Y\mu_j^Y,\,
q_{ji}^Y-\mu_j^Y\mu_i^Y,
q_{ij}^{YZ}-\mu_i^Y\mu_j^Z,\,
q_{ji}^{YZ}-\mu_j^Y\mu_i^Z
\bigr).
\]
Consequently, $\smash{\hat{\bs\rho}_{ij,N}=G(\hat{\bs\eta}_{ij,N})}$ and
$\smash{\bs\rho_{ij}=G(\bs\eta_{ij})}$.
The estimator of the pairwise parameter vector is thus given by
\[
\hat{\bs\theta}_{ij,N}
:=
F\!\big(\hat{\bs\rho}_{ij,N}\big)
=
F\!\big(G(\hat{\bs\eta}_{ij,N})\big).
\]
Thus, the estimator is obtained by first mapping the empirical
cross-products and means into the empirical observable moment vector
through $G$, and then applying the inverse moment map $F$. Explicitly,
\[
\hat m_{ji,N}
=
\frac{\hat\rho_{ij,2,N}^{YZ}}
{\hat\rho_{ij,N}^Y},
\qquad
\hat m_{ij,N}
=
\frac{\hat\rho_{ji,2,N}^{YZ}}
{\hat\rho_{ji,N}^Y},
\]
and
\[
\hat\lambda_{ij,N}
=
\frac{
\hat\rho_{ij,N}^Y\hat\rho_{ji,N}^Y
}{
\hat\rho_{ji,2,N}^{YZ}
}
-
\hat\rho_{ji,N}^Y,
\qquad
\hat\lambda_{ji,N}
=
\frac{
\hat\rho_{ij,N}^Y\hat\rho_{ji,N}^Y
}{
\hat\rho_{ij,2,N}^{YZ}
}
-
\hat\rho_{ij,N}^Y;
\]
cf.\ \eqref{eq:m_inverse} and \eqref{eq:lambda_inverse}.
The full estimator $\smash{\hat{\bs\theta}_N}$ is obtained by collecting
these pairwise estimators over all unordered pairs $\smash{\{i,j\}}$.

Importantly, the resulting estimator $\hat{\bs\theta}_N$ is constructed entirely
from the aggregate observations $(Y_{i,n},Z_{i,n})$. The latent
processes are used only to derive the population moment identities and
to establish the asymptotic properties of the estimator; neither their
observation nor their reconstruction is required for estimation.

\begin{remark}\label{R10}
{\em
The explicit inverse construction above gives a convenient square
system: for each unordered pair $\{i,j\}$, the four moments
$\smash{\rho_{ij}^Y}$, $\smash{\rho_{ji}^Y}$,
$\smash{\rho_{ij,2}^{YZ}}$, and $\smash{\rho_{ji,2}^{YZ}}$ determine the
four parameters $\smash{\lambda_{ij}}$, $\smash{\lambda_{ji}}$,
$\smash{m_{ij}}$, and $\smash{m_{ji}}$. The choice of these moments is
not unique, however. One may, for example, replace some of the lag-2
equations by equations involving the first moments
$\smash{\mu_i^Y}$ and $\smash{\mu_i^Z}$.

It is important in doing so that the resulting equations provide
independent restrictions. A mere count of equations is not sufficient
for identification. Indeed, under the Poisson--Bernoulli specification,
\eqref{eq:BinPo} implies that, for each $i\neq j$,
$\smash{\rho_{ij,2}^{YZ}/\rho_{ij}^Y=m_{ji}}$. Together with the
lag-1 relation $\smash{\rho_{ji}^Y=\mu_{ji}m_{ji}}$, this also determines
$\smash{\mu_{ji}}$. Hence, once sufficiently many lag-2 equations have
been selected, some of the first-moment equations may become implied by
the second-order moment equations and the remaining lag-1 relations.

For example, suppose that $d=4$ and that we use the
$\smash{2d-1=7}$ independent first-moment equations together with the
$\smash{d^2-d=12}$ lag-1 equations based on
$\smash{\rho_{ij}^Y}$, leaving five further equations to be supplied by
lag-2 moments. If the selected lag-2 moments include
$\smash{\rho_{12,2}^{YZ}}$, $\smash{\rho_{13,2}^{YZ}}$, and
$\smash{\rho_{14,2}^{YZ}}$, then the corresponding ratios with
$\smash{\rho_{12}^Y}$, $\smash{\rho_{13}^Y}$, and
$\smash{\rho_{14}^Y}$ determine
$\smash{m_{21}}$, $\smash{m_{31}}$, and $\smash{m_{41}}$.
The corresponding lag-1 relations then determine
$\smash{\mu_{21}}$, $\smash{\mu_{31}}$, and $\smash{\mu_{41}}$.
Consequently, the column-sum relation
$\smash{\sum_{i\neq1}\mu_{i1}=\mu_1^Z}$, which is one of the
first-moment equations, may already be implied by the selected
second-order equations and the remaining moment relations.

Thus, when constructing an estimating system, the relevant requirement
is not merely that the number of selected moment equations matches or
exceeds the number of unknown parameters, but that the selected
equations provide sufficiently many independent restrictions. In
particular, the selection of the lag-2 moments should be made so that
the resulting population moment map has full column rank at the
parameter value under consideration. The explicit pairwise inverse
developed above avoids this selection issue by using a square system
whose inverse is available in closed form.
}
\hfill$\Diamond$
\end{remark}

\begin{remark}\label{R11}
{\em
The inverse relations also reveal a simple relationship between the immigration parameters and the stationary moments. From the stationary mean equation and
\eqref{eq:BinPo}, we have
$\smash{\mu_{ij}=\lambda_{ij}+\rho_{ji}^Y}$ and hence
$\smash{\lambda_{ij}=\mu_{ij}-\rho_{ji}^Y}$. Summing over $j\neq i$
gives
\[
\sum_{j\neq i}\lambda_{ij}
=
\mu_i^Y-\sum_{j\neq i}\rho_{ji}^Y.
\]
Thus, the aggregate first moments together with the lag-1 covariances
identify the total immigration rate associated with each individual.
The individual pairwise immigration parameters are then recovered by
the pairwise inverse map above.
}
\hfill$\Diamond$
\end{remark}
\subsection{Consistency and asymptotic normality}

We next establish the asymptotic behavior of the estimator. The key
probabilistic ingredient is geometric mixing of each latent pair
process.

\begin{lemma}\label{lem:geometric_ergodicity}
Fix $i\neq j$. Suppose that $\lambda_{ij}>0$, $\lambda_{ji}>0$, and
$m_{ij}m_{ji}<1$. Under the Poisson--Bernoulli specification, the
stationary Markov chain $V_n:=(X_{ij,n},X_{ji,n})$ is geometrically
ergodic. In particular, its stationary version is geometrically
$\beta$-mixing and hence geometrically $\alpha$-mixing.
\end{lemma}

\begin{proof}
Consider the two-step skeleton $\widetilde V_n:=V_{2n}$. Let
\[
W(x_1,x_2):=1+x_1+x_2.
\]
Under the Poisson--Bernoulli specification, over two time steps the expected
number of descendants of either type generated by a single individual of that
same type is $m_{ij}m_{ji}$. Hence, since $m_{ij}m_{ji}<1$, there exist
constants $r<1$ and $b<\infty$ such that
\[
\mathbb{E}\!\left[
W(\widetilde V_{n+1})
\mid
\widetilde V_n=(x_1,x_2)
\right]
\leq rW(x_1,x_2)+b.
\]
Thus, the two-step skeleton satisfies a geometric drift condition.

Because the state space is countable and the strictly positive Poisson
immigration rates imply irreducibility and aperiodicity, the finite sublevel
sets of $W$ are small. The Foster--Lyapunov criterion therefore yields
geometric ergodicity of the two-step skeleton; see
\cite[Thm.~15.0.1]{MeynTweedie1993}. Geometric ergodicity of the original
chain then follows. Consequently, the stationary version is geometrically
$\beta$-mixing, and hence geometrically $\alpha$-mixing; see, for example,
\cite[Ch.~16]{MeynTweedie1993}.
\end{proof}

\begin{theorem}[Consistency and asymptotic normality]
Assume that, for every $i\neq j$,
$\lambda_{ij}>0$, $m_{ij}>0$, and $m_{ij}m_{ji}<1$.
Then
\[
\hat{\bs\theta}_N\xrightarrow{\mathbb P}\bs\theta
\qquad\text{and}\qquad
\sqrt{N}\bigl(\hat{\bs\theta}_N-\bs\theta\bigr)
\xrightarrow{d}
\mathcal N(0,\Sigma)
\]
for a finite covariance matrix $\Sigma$.
\end{theorem}

\begin{proof}
For each unordered pair $\{i,j\}$, let
$V_n:=(X_{ij,n},X_{ji,n})$. By
Lemma~\ref{lem:geometric_ergodicity}, there exist constants
$C_{ij}<\infty$ and $r_{ij}\in(0,1)$ such that the corresponding
strong-mixing coefficients satisfy
$\alpha_{ij}(n)\leqslant C_{ij}r_{ij}^n$. Since there are finitely many
unordered pairs, there exist {\it uniform} constants $C<\infty$ and
$r\in(0,1)$ such that
$\alpha_{ij}(n)\leqslant Cr^n$ uniformly over all unordered pairs
$\{i,j\}$.

Under the Poisson--Bernoulli specification, the stationary variables
$X_{ij,n}$ have moments of every order. Hence, the stationary sequences
formed from the finite products of the observed processes that enter
the empirical means and cross-product averages have finite
$(2+\delta)$-th moments for every $\delta>0$. Moreover,
\[
\alpha_{ij}(n)^{\delta/(2+\delta)}
\leqslant
C^{\delta/(2+\delta)}
r^{n\delta/(2+\delta)},
\]
and therefore the required summability condition on the strong-mixing
coefficients is satisfied. Thus, the moment and mixing conditions of the
central limit theorem for strongly mixing sequences in
\cite[Thm.\ 18.5.3]{ibragimov1971} are fulfilled for the stationary
sequences whose sample averages form the components of
$\hat{\bs\eta}_{ij,N}$. Consequently, for each unordered pair
$\{i,j\}$, the empirical means and cross-product averages collected in
$\hat{\bs\eta}_{ij,N}$ satisfy a joint central limit theorem.

For each unordered pair $\{i,j\}$, let $\bs\eta_{ij}$ and
$\hat{\bs\eta}_{ij,N}$ be as defined above. Stacking these vectors over
all unordered pairs gives the global vectors $\hat{\bs\eta}_N$ and
$\bs\eta$. Since there are $\binom{d}{2}$ unordered pairs and each pair
contributes eight components, we have that 
$\hat{\bs\eta}_N$ and $\bs\eta$ have dimension ${8\binom{d}{2}}
={4d(d-1)}.$ The pairwise joint central limit theorems, together with the
independence of processes corresponding to distinct unordered pairs,
and the Cram\'er--Wold device therefore give
\[
\sqrt{N}\bigl(\hat{\bs\eta}_N-\bs\eta\bigr)
\xrightarrow{d}
\mathcal N(0,\Sigma^\eta).
\]
The independence across unordered pairs implies that
$\Sigma^\eta$ is obtained by combining the corresponding pairwise
asymptotic covariance matrices. The
finiteness of its entries follows by Lemma \ref{lem:finite_Sigma_eta}.
The map $G:\mathbb R^8\to\mathbb R^4$ introduced in the previous
subsection acts separately on each unordered pair. Its blockwise
extension to the full collection of unordered pairs will be denoted by
\[
\mathcal G:
\mathbb R^{4d(d-1)}
\longrightarrow
\mathbb R^{2d(d-1)}.
\]
Thus, $\smash{\hat{\bs\rho}_N=\mathcal G(\hat{\bs\eta}_N)}$ and
$\smash{\bs\rho=\mathcal G(\bs\eta)}$.
Since $G$ is continuously differentiable, so is $\mathcal G$, with
Jacobian $J_{\mathcal G}(\bs\eta)$ given by the block-diagonal matrix
whose blocks are the Jacobian $J_G(\bs\eta_{ij})$ of the pairwise map.
In particular,
\[
J_{\mathcal G}(\bs\eta)\in
\mathbb R^{2d(d-1)\times 4d(d-1)}.
\]
The multivariate delta method therefore gives
\[
\sqrt{N}\bigl(\hat{\bs\rho}_N-\bs\rho\bigr)
\xrightarrow{d}
\mathcal N(0,\Sigma^\rho),
\qquad
\Sigma^\rho
=
J_{\mathcal G}(\bs\eta)\,
\Sigma^\eta\,
J_{\mathcal G}(\bs\eta)^\top,
\]
where $\smash{\Sigma^\eta\in\mathbb R^{4d(d-1)\times4d(d-1)}}$ and
$\smash{\Sigma^\rho\in\mathbb R^{2d(d-1)\times2d(d-1)}}$.

By the assumptions $\lambda_{ij}>0$ and $m_{ij}>0$, the population
moments appearing in the denominators of the inverse map are strictly
positive. Hence the pairwise inverse map $F$ introduced above is
continuously differentiable in a neighborhood of each
$\bs\rho_{ij}$. Let $\mathcal F$ denote its blockwise extension to the
full collection of unordered pairs,
\[
\mathcal F:
\mathbb R^{2d(d-1)}
\longrightarrow
\mathbb R^{2d(d-1)}.
\]
Thus, $\smash{\hat{\bs\theta}_N=\mathcal F(\hat{\bs\rho}_N)}$ and
$\smash{\bs\theta=\mathcal F(\bs\rho)}$.
Since $\mathcal F$ is continuously differentiable, its Jacobian satisfies
\[
J_{\mathcal F}(\bs\rho)\in
\mathbb R^{2d(d-1)\times 2d(d-1)}.
\]
A second application of the multivariate delta method gives
\[
\sqrt{N}\bigl(\hat{\bs\theta}_N-\bs\theta\bigr)
\xrightarrow{d}
\mathcal N(0,\Sigma),
\qquad
\Sigma
=
J_{\mathcal F}(\bs\rho)
\,\Sigma^\rho\,
J_{\mathcal F}(\bs\rho)^\top,
\]
where $\smash{\Sigma\in\mathbb R^{2d(d-1)\times2d(d-1)}}$.
Equivalently, since the estimator is obtained by the composition
$\mathcal F\circ\mathcal G$, the asymptotic covariance can be written as
\[
\Sigma
=
J_{\mathcal F}(\bs\rho)
J_{\mathcal G}(\bs\eta)\,
\Sigma^\eta\,
J_{\mathcal G}(\bs\eta)^\top
J_{\mathcal F}(\bs\rho)^\top.
\]

Finally, the law of large numbers gives consistency of
$\smash{\hat{\bs\eta}_N}$. Since $\mathcal G$ is continuous, this implies
consistency of $\smash{\hat{\bs\rho}_N=\mathcal G(\hat{\bs\eta}_N)}$. The
continuity of $\mathcal F$ at $\smash{\bs\rho}$ then gives consistency of
$\smash{\hat{\bs\theta}_N=\mathcal F(\hat{\bs\rho}_N)}$. This completes the
proof.
\end{proof}

\begin{remark}\label{R12}
{\em
The proof separates the asymptotic analysis into two parts. The
probabilistic part establishes a joint central limit theorem for the
empirical quantities collected in $\smash{\hat{\bs\eta}_N}$. The
remaining steps are deterministic: the map $\mathcal G$ converts these
quantities into the empirical moment vector, while the inverse map
$\mathcal F$ converts the latter into the parameter estimator. Thus,
once the joint central limit theorem for $\smash{\hat{\bs\eta}_N}$ is
available, consistency and asymptotic normality follow from the
continuity and differentiability of $\mathcal G$ and $\mathcal F$.
}
\hfill$\Diamond$
\end{remark}

\section{Numerical Experiments}\label{sec:num}

In this section illustrate the performance of the proposed estimation procedure through numerical experiments. We begin with a setting satisfying the assumptions of the main asymptotic result,
\begin{equation}\label{eq:assu}
\lambda_{ij}>0,\qquad
m_{ij}>0,\qquad
m_{ij}m_{ji}<1,
\qquad i\neq j.
\end{equation}
We include two examples that illustrate the consistency and $N^{-1/2}$ scaling established theoretically above.

We then investigate the performance of the procedure beyond the formal scope of the theorem. Although the positivity conditions on $\lambda_{ij}$ and $m_{ij}$ are required for the theoretical results, they are not intrinsic limitations of the estimation procedure. We therefore consider an additional example in which some of these conditions are violated.

\subsection{Numerical demonstration of the proposed approach}

We conduct two experiments under the assumptions \eqref{eq:assu} stated above. Experiment~1 evaluates the accuracy of the parameter estimators, while Experiment~2 examines how estimation accuracy varies with the number of observations.

In both experiments, we consider a network consisting of $d=5$
individuals, with system parameters given by
\[
\bs{\lambda} =
\begin{pmatrix}
0.000 & 0.703 & 0.823 & 0.711 & 0.679 \\
0.688 & 0.000 & 1.042 & 0.899 & 0.860 \\
0.579 & 0.748 & 0.000 & 0.756 & 0.723 \\
0.715 & 0.925 & 1.083 & 0.000 & 0.894 \\
0.732 & 0.946 & 1.108 & 0.956 & 0.000
\end{pmatrix},
\]
and
\[
\bs{m} =
\begin{pmatrix}
0.000 & 0.402 & 0.470 & 0.406 & 0.388 \\
0.393 & 0.000 & 0.595 & 0.514 & 0.491 \\
0.331 & 0.428 & 0.000 & 0.432 & 0.413 \\
0.409 & 0.528 & 0.619 & 0.000 & 0.511 \\
0.418 & 0.541 & 0.633 & 0.547 & 0.000
\end{pmatrix}.
\]

\medskip

\noindent\textbf{\textit{Experiment 1.}}
The main objective of this experiment is to assess the accuracy of the
proposed parameter estimators. We perform $R=1\,000$ independent
simulation runs, each consisting of $N=250\,000$ observation times. For
each pair $i\neq j$, we compute the sample mean and sample standard
deviation from the resulting $1\,000$ estimates of $\hat m_{ij}$ and
$\hat\lambda_{ij}$.

\medskip

\begin{enumerate}
    \item The sample mean and sample standard deviation of
    $\hat{\bs{\lambda}}$ are
    \begin{align*}
    {\bs{\tilde{m}}_{\bs{\lambda}}^{\mathrm{avg}}}
    &=
    \begin{pmatrix}
    0.000 & 0.706 & 0.833 & 0.713 & 0.686 \\
    0.692 & 0.000 & 1.045 & 0.900 & 0.861 \\
    0.582 & 0.747 & 0.000 & 0.757 & 0.722 \\
    0.719 & 0.926 & 1.083 & 0.000 & 0.893 \\
    0.738 & 0.948 & 1.111 & 0.957 & 0.000
    \end{pmatrix},
    \\[2mm]
    \bs{s}_{\bs{\lambda}}
    &=
    10^{-1}
    \begin{pmatrix}
    0.00 & 0.72 & 0.76 & 0.66 & 0.73 \\
    0.71 & 0.00 & 0.50 & 0.49 & 0.53 \\
    0.72 & 0.58 & 0.00 & 0.56 & 0.57 \\
    0.68 & 0.51 & 0.48 & 0.00 & 0.49 \\
    0.67 & 0.50 & 0.49 & 0.49 & 0.00
    \end{pmatrix}.
    \end{align*}

    \item The sample mean and sample standard deviation of
    $\hat{\bs{m}}$ are
    \begin{align*}
    {\bs{\tilde{m}}_{\bs{m}}^{\mathrm{avg}}}
    &=
    \begin{pmatrix}
    0.000 & 0.402 & 0.468 & 0.406 & 0.387 \\
    0.392 & 0.000 & 0.595 & 0.514 & 0.491 \\
    0.331 & 0.428 & 0.000 & 0.432 & 0.414 \\
    0.409 & 0.528 & 0.619 & 0.000 & 0.511 \\
    0.416 & 0.541 & 0.633 & 0.547 & 0.000
    \end{pmatrix},
    \\[2mm]
    \bs{\nu}_{\bs{m}}
    &=
    10^{-1}
    \begin{pmatrix}
    0.00 & 0.23 & 0.28 & 0.21 & 0.22 \\
    0.22 & 0.00 & 0.15 & 0.11 & 0.12 \\
    0.18 & 0.11 & 0.00 & 0.11 & 0.10 \\
    0.21 & 0.12 & 0.14 & 0.00 & 0.11 \\
    0.22 & 0.13 & 0.15 & 0.12 & 0.00
    \end{pmatrix}.
    \end{align*}
\end{enumerate}

The simulation results indicate that the estimators recover the
underlying parameter values accurately. In particular, the sample means
are close to the corresponding true parameter values, while the
standard deviations are relatively small. These results provide a
numerical illustration of the finite-sample performance of the proposed
estimation procedure for the network considered here.

We next examine the empirical sampling distributions of the estimators.
Since the total number of parameters is as large as 
$
2d^2-2d=40
$
in this example, we restrict attention to the representative parameter
pair $(\lambda_{15},m_{15})$. The remaining parameter pairs exhibit
similar qualitative behavior.

\begin{figure}[!htb]
    \centering
    \includegraphics[width=0.7\textwidth]{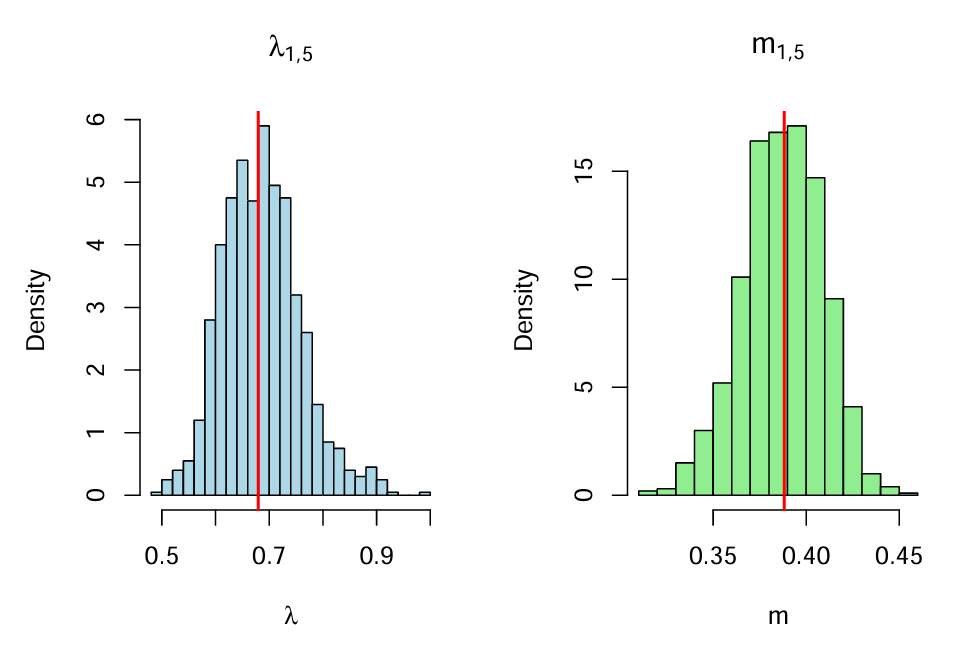}
    \caption{Histograms of the estimators of $\lambda_{15}=0.679$ and
    $m_{15}=0.388$, based on $R=1\,000$ independent experiments, each
    with $N=250\,000$ observations. The red vertical line in each
    histogram indicates the corresponding true parameter value.}
    \label{fig:hist_lambda}
\end{figure}

\begin{figure}[H]
    \centering
    \includegraphics[width=0.7\linewidth]{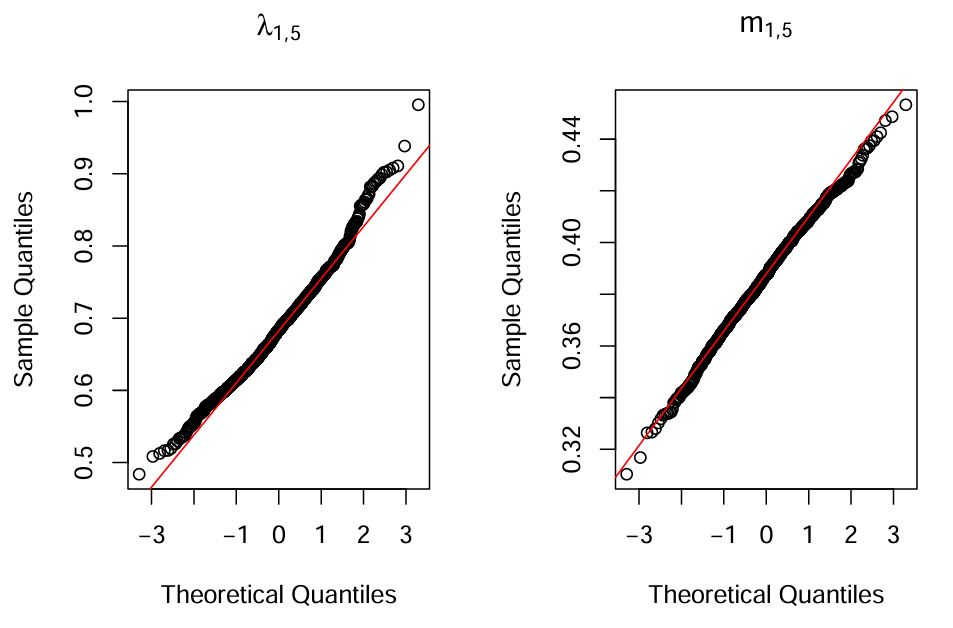}
    \caption{QQ-plots of the estimators of $\lambda_{15}$ and $m_{15}$
    based on $R=1\,000$ independent experiments, each with
    $N=250\,000$ observations.}
    \label{fig:qq_lambda15}
\end{figure}

Figure~\ref{fig:hist_lambda} displays the corresponding histograms. The
estimates are concentrated around their respective true parameter
values, illustrating the accuracy of the estimators. Moreover, the
histograms suggest that the sampling distributions are approximately
normal. This observation is further illustrated by the QQ-plots shown
in Figure~\ref{fig:qq_lambda15}.

\medskip

\noindent\textbf{\textit{Experiment 2.}}
The primary objective of this experiment is to investigate how the
number of observations affects the accuracy of the estimators. We
consider six sample sizes,
\[
N\in
\left\{
10^5,\;
250\,000,\;
500\,000,\;
750\,000,\;
10^6,\;
2\times10^6
\right\},
\]
and compute the standard deviation of each parameter estimator for
every value of $N$.

\begin{figure}[H]
    \centering
    \includegraphics[width=1\linewidth]{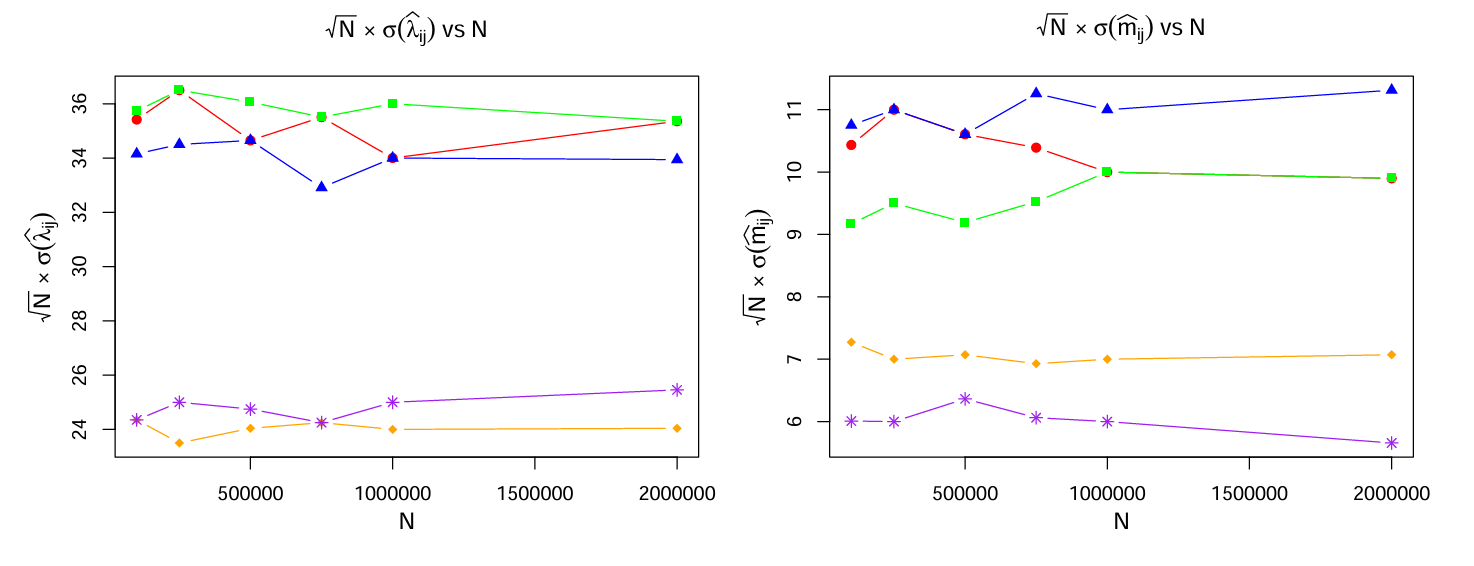}
    \caption[Scaled standard deviations for five representative parameter pairs]{
Scaled standard deviations of the estimators
$\sqrt{N}\,\sigma(\widehat{\lambda}_{ij})$ (left) and
$\sqrt{N}\,\sigma(\widehat{m}_{ij})$ (right), as functions
of the sample size $N$ for five representative parameter pairs.
The sample sizes are $N=10^5$, $250\,000$, $500\,000$,
$750\,000$, $10^6$, and $2\times10^6$, with $R=1\,000$
independent experiments for each sample size.
The five representative parameter pairs are indicated by the
following markers:
\protect\Rred~$(1,5)$;
\protect\Rblue~$(5,1)$;
\protect\Rgreen~$(3,1)$;
\protect\Rorange~$(5,3)$; and
\protect\Rpurple~$(4,2)$.
}
    \label{fig:Nvalues}
\end{figure}

For ease of presentation, we select five representative parameter pairs
and plot the corresponding standard deviations as functions of $N$.
The results for both $\hat\lambda_{ij}$ and $\hat m_{ij}$ are displayed
in Figure~\ref{fig:Nvalues}, where
${\sigma}(
\hat{\lambda}_{ij})$ and ${\sigma}(\hat{m}_{ij})$ denote the standard deviations of
the corresponding estimators.

As the sample size increases, the standard deviations decrease. In
addition, the values of
${\sigma}(\hat{\lambda}_{ij})\sqrt{N}$ and ${\sigma}(\hat{m}_{ij})\sqrt{N}$ shown in
Figure~\ref{fig:Nvalues} remain approximately constant as \(N\) increases, consistent with the $N^{-1/2}$ scaling
predicted by the asymptotic theory. A similar pattern is observed for
the other parameter pairs.

\subsection{General network structures}\label{ssec:identif}

In the previous analysis, we assumed that all parameters are strictly
positive. Under this assumption, the lag-1 moments
$\smash{\rho_{ij}^Y}$ and the lag-2 moments
$\smash{\rho_{ij,2}^{YZ}}$ are all strictly positive. We next consider
cases in which some of the parameters are zero.

For a given pair $\{i,j\}$, the four associated parameters are
$\smash{\lambda_{ij}}$, $\smash{\lambda_{ji}}$,
$\smash{m_{ij}}$, and $\smash{m_{ji}}$. The corresponding four observable
moments are $\smash{\rho_{ij}^Y}$, $\smash{\rho_{ji}^Y}$,
$\smash{\rho_{ij,2}^{YZ}}$, and $\smash{\rho_{ji,2}^{YZ}}$. The
relationships in \eqref{eq:BinPo} show that the vanishing pattern of the
lag-1 moments provides useful information about the underlying parameter
configuration. We consider the relevant cases in turn.

\medskip

\noindent\textbf{--- {\it All four moments are positive.}}
Suppose that
$\smash{\rho_{ij}^Y>0}$, $\smash{\rho_{ji}^Y>0}$,
$\smash{\rho_{ij,2}^{YZ}>0}$, and
$\smash{\rho_{ji,2}^{YZ}>0}$. This can occur when all four parameters
are strictly positive, as in the setting considered in the previous
subsection, or when one of the immigration parameters is zero.

For example, suppose, without loss of generality, that
$\smash{\lambda_{ij}=0}$. Since $\smash{\rho_{ij}^Y>0}$ and
$\smash{\rho_{ji}^Y>0}$, the stationary means are positive, and hence
$\smash{m_{ij}>0}$ and $\smash{m_{ji}>0}$. Although there is no
exogenous immigration from $i$ to $j$, communication from $j$ to $i$
can generate a response from $i$ to $j$. Thus, the lagged moments need
not vanish, and all four moments may remain positive. In this case, the
inverse formulas derived above continue to apply and recover
$\smash{\lambda_{ij}=0}$ together with the remaining three parameters.

\medskip

\noindent\textbf{--- {\it All four moments are zero.}}
Next, consider the case
$\smash{\rho_{ij}^Y=\rho_{ji}^Y=\rho_{ij,2}^{YZ}
=\rho_{ji,2}^{YZ}=0}$. In this case, the moment conditions alone do not
uniquely determine the underlying parameter configuration. For example,
one possibility is that there is no exogenous interaction between $i$
and $j$, corresponding to
$\smash{\lambda_{ij}=\lambda_{ji}=0}$. Another possibility is that both
response parameters vanish, corresponding to
$\smash{m_{ij}=m_{ji}=0}$. These distinct parameter configurations can
therefore generate the same collection of zero-valued moments. {Note that these are just two illustrative examples; other parameter configurations may also lead to the same zero-valued moments.}
Consequently, the network structure is not identifiable from these
moments alone.

\medskip

\noindent\textbf{--- {\it Exactly one lag-1 moment {and two lag-2 moments} are zero.}}
Finally, suppose, without loss of generality, that
$\smash{\rho_{ij}^Y=0}$ and $\smash{\rho_{ji}^Y>0}$. Since
$\smash{\rho_{ji}^Y>0}$ implies
$\smash{\mu_{ji}>0}$ and $\smash{m_{ji}>0}$, the stationary mean
equation implies $\smash{\mu_{ij}>0}$. {Therefore, in this case, since \(\mu_{ij}>0\), Equation \eqref{eq:BinPo} implies}
\[
\rho_{ij}^Y=0
\quad\Longleftrightarrow\quad
m_{ij}=0.
\]
Thus, the zero value of
$\smash{m_{ij}}$ is directly identified from the vanishing lag-1 moment.

The remaining parameters associated with the pair cannot, however, be
recovered from the inverse formulas based on \eqref{eq:BinPo}, {since the relevant ratios involve division by $\rho_{ji}^Y$ and the lag-2 moments, which are all zero in this case.}
Identification
must therefore be supplemented by the stationary mean equations and
additional identifying restrictions. In particular, when several pairs
have a zero response parameter, the corresponding system of mean
equations must contain sufficiently many independent restrictions to
recover the remaining immigration parameters. A convenient sufficient
condition is that the relevant collection of pairs form an acyclic set,
together with the requirement that the number of remaining unknown
immigration parameters not exceed the number of available independent
mean equations.

\medskip

In summary, zero-valued parameters can lead to non-identification in
two distinct ways. When all four observable moments vanish, different
underlying parameter configurations may produce the same observable
moment vector. When a response parameter is zero, the explicit inverse
map developed above breaks down because the corresponding lag-1 moment
vanishes, and additional mean restrictions are required. By contrast, a
zero immigration parameter does not by itself create an identification
problem: as long as the relevant response parameters and stationary
means remain positive, the pairwise inverse map continues to recover
the zero immigration parameter. The following numerical experiment illustrates the estimation procedure for a sparse network in which the four moments $\smash{\rho_{ij}^Y}$, $\smash{\rho_{ji}^Y}$, $\smash{\rho_{ij,2}^{YZ}}$, and $\smash{\rho_{ji,2}^{YZ}}$ are all positive, demonstrating that the model is identifiable.

\medskip


We next present an additional simulation experiment to illustrate the identification of a sparse network. Before presenting the results, we clarify the notations used in the network diagrams. Dashed edges represent the $m-$parameters, while solid edges correspond to $\lambda-$parameters. Arrowheads indicate the direction of the corresponding parameters. Solid edges with a single arrowhead are shown in brown to highlight the presence of a zero-valued $\lambda$-parameter.
\begin{figure}[H]
    \centering
    \includegraphics[width=0.3\textwidth]{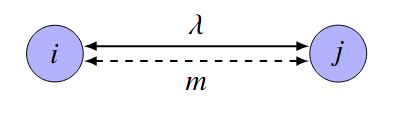}
\end{figure}


\noindent\noindent\textbf{\textit{Experiment 3.}} The main aim of this experiment is to illustrate a scenario in which sparse network can be recovered and to evaluate the reliability of the estimators, particularly for parameters whose true values are zero. As in the previous experiments, we consider a network with $d=5$ nodes and perform $R=1\,000 $ independent simulation runs, each consisting of $N=250\,000$ observation times. The system parameters to be estimated are listed below:
\medskip

\begin{minipage}[c]{0.5\textwidth}
\[
\begin{aligned}
\bs{\lambda} &=
\begin{pmatrix}
0.000 & 0.703 & 0.823 & 0.711 & 0.679 \\
0.688 & 0.000 & 0.000 & 0.000 & 0.000 \\
0.579 & 0.748 & 0.000 & 0.756 & 0.723 \\
0.000 & 0.925 & 0.000 & 0.000 & 0.000 \\
0.732 & 0.946 & 0.000 & 0.956 & 0.000
\end{pmatrix} \\[1em]
\bs{m} &=
\begin{pmatrix}
0.000 & 0.402 & 0.470 & 0.406 & 0.388 \\
0.393 & 0.000 & 0.595 & 0.514 & 0.491 \\
0.331 & 0.428 & 0.000 & 0.432 & 0.413 \\
0.409 & 0.528 & 0.619 & 0.000 & 0.511 \\
0.418 & 0.541 & 0.633 & 0.547 & 0.000
\end{pmatrix}.
\end{aligned}
\]
\end{minipage}
\hfill
\begin{minipage}[c]{0.4\textwidth}
\centering
\begin{tikzpicture}[>=latex]
    \fivenodes
    \twoway{1}{2}
    \twoway{1}{3}
    \twoway{1}{5}
    \oneway{3}{2}
    \oneway{4}{2}
    \oneway{5}{2}
    \oneway{1}{4}
    \oneway{3}{4}
    \oneway{5}{4}
    \oneway{3}{5}
\end{tikzpicture}
\end{minipage}

\medskip

We first report the empirical means and standard deviations of the estimators across the $1 \,000$ independent replications, to assess the accuracy of the estimates of both the $\lambda$- and $m$-parameters. 

\medskip

\begin{enumerate} 
\item The sample mean and sample standard deviation of $\hat{\bs{\lambda}}$ are
\begin{align*}
{\bs{\tilde{m}}_{\bs{\lambda}}^{\mathrm{avg}}}
&=
\begin{pmatrix}
0.000 & 0.704 & 0.826 & 0.724 & 0.679 \\
0.692 & 0.000 & 0.001 & 0.000 & 0.000 \\
0.580 & 0.750 & 0.000 & 0.758 & 0.728 \\
0.001 & 0.925 & 0.000 & 0.000 & 0.000 \\
0.732 & 0.945 & 0.000 & 0.959 & 0.000
\end{pmatrix},
\\[2mm]
\bs{s}_{\bs{\lambda}}
&=
10^{-1}
\begin{pmatrix}
0.00 & 0.42 & 0.51 & 0.96 & 0.50 \\
0.57 & 0.00 & 0.18 & 0.16 & 0.22 \\
0.40 & 0.46 & 0.00 & 0.43 & 0.57 \\
0.25 & 0.43 & 0.16 & 0.00 & 0.19 \\
0.45 & 0.41 & 0.16 & 0.37 & 0.00
\end{pmatrix}.
\end{align*}

\medskip

\item The sample mean and sample standard deviation of $\hat{\bs{m}}$ are
\begin{align*}
{\bs{\tilde{m}}_{\bs{m}}^{\mathrm{avg}}}
&=
\begin{pmatrix}
0.000 & 0.402 & 0.470 & 0.405 & 0.389 \\
0.393 & 0.000 & 0.595 & 0.514 & 0.491 \\
0.330 & 0.427 & 0.000 & 0.432 & 0.412 \\
0.408 & 0.529 & 0.618 & 0.000 & 0.510 \\
0.418 & 0.542 & 0.633 & 0.546 & 0.000
\end{pmatrix},
\\[2mm]
\bs{s}_{\bs{m}}
&=
10^{-1}
\begin{pmatrix}
0.00 & 0.13 & 0.19 & 0.43 & 0.15 \\
0.18 & 0.00 & 0.19 & 0.10 & 0.12 \\
0.10 & 0.18 & 0.00 & 0.16 & 0.22 \\
0.23 & 0.17 & 0.16 & 0.00 & 0.10 \\
0.15 & 0.16 & 0.16 & 0.14 & 0.00
\end{pmatrix}.
\end{align*}
\end{enumerate}
\medskip
The results indicate that the estimators exhibit low sampling variability, demonstrating the stability of the proposed estimation algorithm.

\begin{figure}
    \centering
    \includegraphics[width=0.7\textwidth]{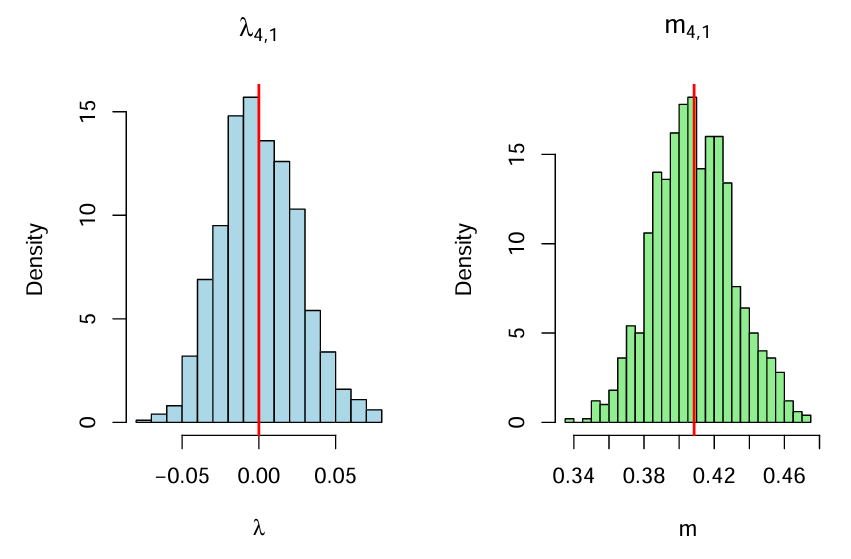}
    \caption{Histograms pertaining to estimates of $\lambda_{41}=0$ and $m_{41}=0.409$, based on $R=1000$ experiments, each with $N=250\,000$ observations. In each histogram, the red vertical line indicates the true parameter value.}
    \label{fig:hist_lambda_zero}
\end{figure}

\begin{figure}
    \centering
    \includegraphics[width=0.7\linewidth]{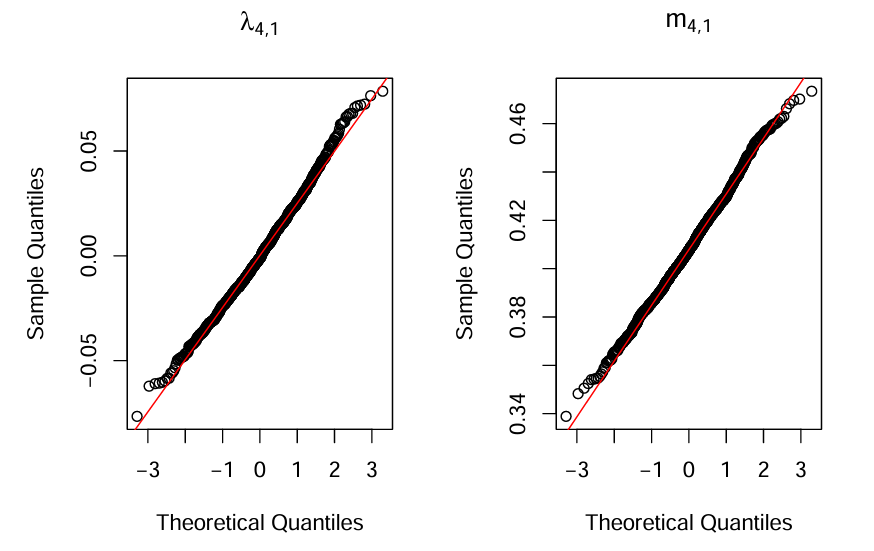}
    \caption{QQ-plots of the estimators of $\lambda_{41}$ and $m_{41}$ based on $R=1000$ experiments, each with $N=250\,000$ observations.}
    \label{fig:qq_lambda41}
\end{figure}

To further illustrate the behavior of the estimators, we present the histogram and QQ-plot for the representative pair of parameters $(\lambda_{41}, \, m_{41})$ in Figures~\ref{fig:hist_lambda_zero} and~\ref{fig:qq_lambda41}. This pair provides a typical example of the estimation procedure in the presence of zero-valued parameters. The remaining identifiable estimators exhibit the same qualitative behavior.

The histograms suggest that the empirical sampling distributions of the estimators are approximately normal and centered around the true parameter values, indicating good estimation accuracy. The corresponding Q–Q plots provide further evidence of approximate normality.
\medskip

\section{Concluding Remarks}\label{sec: conc}

In this paper, we studied parameter estimation for a stochastic
interaction model when only aggregate observations are available and the
underlying pairwise interaction processes are unobserved. By exploiting
the stationary moment structure of the model, we showed that the latent
parameters can be identified and estimated directly from aggregate data,
without reconstructing the latent processes or maximizing a
high-dimensional likelihood. Under the Poisson--Bernoulli specification,
we obtained an explicit pairwise inverse map from observable moments to
the model parameters and established consistency and asymptotic
normality of the resulting estimator. Numerical experiments illustrated
the finite-sample performance of the procedure and its expected
$N^{-1/2}$ scaling.

The identification analysis also clarifies the role of zero-valued
parameters. In particular, a zero immigration parameter does not
necessarily prevent identification, whereas configurations involving
zero response parameters or vanishing observable moments may require
additional restrictions. Thus, identification depends not only on the
number of available moment equations, but also on the structure and
independence of the selected moments.

Although communication networks provide the motivating example, the
framework applies more broadly to partially observed interacting count
processes in which latent pairwise interactions generate observable
aggregate quantities. For example, aggregate activity may be observable
in online social networks even when detailed interaction data are
unavailable because of privacy, data-storage, or data-collection
constraints.

Several directions for further research are natural.

\begin{itemize}
    \item[$\circ$] \emph{Sparse networks.}
    The analysis in Section~\ref{ssec:identif} shows that zero-valued
    parameters can create additional identification problems. An
    important direction is therefore to develop estimation methods that
    explicitly exploit sparsity while accounting for these
    identification issues.

    \item[$\circ$] \emph{More general interaction structures.}
    The present model is based on pairwise interactions and specific
    immigration and offspring mechanisms. It would be natural to allow
    for interactions involving more than two individuals, time-varying
    parameters, or more general observation and censoring mechanisms.

    \item[$\circ$] \emph{Alternative estimating systems.}
    The explicit inverse map used here provides an exactly identified
    estimating system. One could instead combine first moments with
    additional lagged moments to obtain overidentified estimating
    systems, which could be useful for specification testing or for
    improving statistical efficiency.

    \item[$\circ$] \emph{Growing networks.}
    Our asymptotic results treat the number of individuals as fixed.
    Extending the theory to settings in which the network size grows
    with the number of observations would be important for applications
    involving large networks.
\end{itemize}

\bibliographystyle{plain}
\bibliography{literature}

\appendix 

\section{Two lemmas}

\begin{lemma}\label{L3}
Fix $i\neq j$. Assume that  $m_{ij}m_{ji}<1$. Define the two interlaced processes
\[
V_{ij,n}
:=
\bigl(X_{ij,2n},X_{ji,2n+1}\bigr),
\qquad
V_{ji,n}
:=
\bigl(X_{ji,2n},X_{ij,2n+1}\bigr),
\qquad n\geqslant 0.
\]
Then the processes $(V_{ij,n})_{n\geqslant 0}$ and
$(V_{ji,n})_{n\geqslant 0}$ are independent.
\end{lemma}
 
\begin{proof}
Define the collections of primitive variables
\begin{align*}
\mathscr U_{ij}
&:=
\{I_{ij,2k},I_{ji,2k-1},
\xi_{ji,2k},\xi_{ij,2k-1}:k\in\mathbb Z\},\\
\mathscr U_{ji}
&:=
\{I_{ji,2k},I_{ij,2k-1},
\xi_{ij,2k},\xi_{ji,2k-1}:k\in\mathbb Z\}.
\end{align*}
These collections are disjoint. Hence, by the independence assumptions
on the primitive immigration and offspring variables,
$\sigma(\mathscr U_{ij})$ and $\sigma(\mathscr U_{ji})$ are independent.

The recursion \eqref{eq:model} shows that $X_{ij,2n}$ and
$X_{ji,2n+1}$ depend only on the primitive variables in
$\mathscr U_{ij}$. Thus $V_{ij,n}$ is measurable with respect to
$\sigma(\mathscr U_{ij})$. Likewise, $V_{ji,n}$ is measurable with
respect to $\sigma(\mathscr U_{ji})$.

In the stationary regime, the condition $m_{ij}m_{ji}<1$ ensures that
the two-step recursion is contractive in mean, so the stationary
process admits a causal representation in terms of the primitive
variables extending to the infinite past. Consequently,
\[
\sigma(V_{ij,n}:n\geqslant  0)
\subseteq\sigma(\mathscr U_{ij}),
\qquad
\sigma(V_{ji,n}:n\geqslant  0)
\subseteq\sigma(\mathscr U_{ji}).
\]
Since the two sigma-fields on the right-hand side are independent, the
two processes $(V_{ij,n})_{n\geqslant 0}$ and
$(V_{ji,n})_{n\geqslant 0}$ are independent.
\end{proof}

\begin{lemma}\label{lem:finite_Sigma_eta}
Assume that, for every $i\neq j$,
$\lambda_{ij}>0$, $m_{ij}>0$, and $m_{ij}m_{ji}<1$.
Under the Poisson--Bernoulli specification, the asymptotic covariance matrix
$\Sigma^\eta$ of $\hat{\bs\eta}_N$ has finite entries.
\end{lemma}

\begin{proof}
For each unordered pair $\{i,j\}$, the components of
$\hat{\bs\eta}_{ij,N}$ are empirical averages of the four products
\[
Y_{i,n}Y_{j,n+1},\quad
Y_{j,n}Y_{i,n+1},\quad
Y_{i,n}Z_{j,n+2},\quad
Y_{j,n}Z_{i,n+2},
\]
and the four individual observations
\[
Y_{i,n},\quad Y_{j,n},\quad Z_{i,n},\quad Z_{j,n}.
\]
Collect the corresponding stationary summands, over all unordered
pairs, into a single vector $\bs H_n$. Each component of $\hat{\bs\eta}_N$ is the average of the corresponding
component of $\bs H_n$, over $n=1,\ldots,N$ for the individual
observations, over $n=1,\ldots,N-1$ for the lag-1 products, and over
$n=1,\ldots,N-2$ for the lag-2 products. Since the differences between
these summation limits and $N$ are fixed, they do not affect the
asymptotic covariance matrix, and all components may be treated as
averages with denominator $N$ for the purpose of the asymptotic
analysis. Observe that the entries of the asymptotic covariance matrix are the corresponding
long-run covariance sums:
\begin{equation}\label{eq:sigma_eta}
\Sigma^\eta_{k\ell}
=
\sum_{n\in\mathbb Z}
\operatorname{Cov}(H_0^{(k)},H_n^{(\ell)}),
\end{equation}
where $H_n^{(k)}$ denotes the $k$-th component of $\bs H_n$.

Under the Poisson--Bernoulli specification, the latent variables have
moments of every order. Since $Y_{i,n}$ and $Z_{i,n}$ are finite sums
of latent variables, they also have moments of every order. By
Hölder's inequality, the same holds for the products
$Y_{i,n}Y_{j,n+1}$ and $Y_{i,n}Z_{j,n+2}$ entering $\bs H_n$.
Consequently, since $\bs H_n$ has finitely many components, we obtain
\[
\mathbb E\|\bs H_0\|^{2+\delta}<\infty
\]
for every $\delta>0$.
By Lemma~\ref{lem:geometric_ergodicity}, the latent pair processes are
geometrically strongly mixing. Their finite joint collection is
therefore geometrically strongly mixing, and so is $\bs H_n$, since
each component of $\bs H_n$ is a function of latent variables at only
finitely many time points. Hence, for some $C<\infty$ and $r\in(0,1)$, we have {geometric mixing}:
\[
\alpha_H(n)\leqslant Cr^n,\qquad k\geqslant1,
\]
where $\alpha_H(n)$ is the strong-mixing coefficient corresponding to the process $(\bs H_n)_{n\in{\mathbb Z}}$.

For any two component indices $k$ and $\ell$, we consider the covariance
between the corresponding components of $\bs H_0$ and $\bs H_n$,
namely,
$ 
\operatorname{Cov}\bigl(H_0^{(k)},H_n^{(\ell)}\bigr).
$
By Davydov's covariance inequality \cite{Davydov1968}, for every $\delta>0$,
\[
\big|\operatorname{Cov}
  (H_0^{(k)},H_n^{(\ell)})\big|
\leqslant C_{k,\ell}\,\alpha_H(n)^{\delta/(2+\delta)}
\]
for a finite constant $C_{k,\ell}$. Geometric mixing consequently implies
\[
\sum_{n\in\mathbb Z}
\big|\operatorname{Cov}
  (H_0^{(k)},H_n^{(\ell)})\big|<\infty
\]
for every pair of components $k,\ell$.
Recall from \eqref{eq:sigma_eta} that the entries of the asymptotic covariance matrix can be written as sums of covariances $ 
\operatorname{Cov}\bigl(H_0^{(k)},H_n^{(\ell)}\bigr).
$
As all these sums converge absolutely, every entry of $\Sigma^\eta$
is finite.
\end{proof}

 \end{document}